\documentclass[11pt,fleqn]{article}
\usepackage{amsfonts, epsfig}
\newcommand{\ol}{\setlength{\itemsep}{0pt.}\begin{enumerate}}
\newcommand{\eol}{\end{enumerate}\setlength{\itemsep}{-\parsep}}
\newcommand{\ignore}[1]{}
\title{Lower bounds for designs in symmetric spaces}
\author{Noa Eidelstein and Alex Samorodnitsky}

\begin{document}
\date{}
\maketitle


\newtheorem{THEOREM}{Theorem}[section]
\newenvironment{theorem}{\begin{THEOREM} \hspace{-.85em} {\bf :}
}%
                        {\end{THEOREM}}
\newtheorem{LEMMA}[THEOREM]{Lemma}
\newenvironment{lemma}{\begin{LEMMA} \hspace{-.85em} {\bf :} }%
                      {\end{LEMMA}}
\newtheorem{COROLLARY}[THEOREM]{Corollary}
\newenvironment{corollary}{\begin{COROLLARY} \hspace{-.85em} {\bf
:} }%
                          {\end{COROLLARY}}
\newtheorem{PROPOSITION}[THEOREM]{Proposition}
\newenvironment{proposition}{\begin{PROPOSITION} \hspace{-.85em}
{\bf :} }%
                            {\end{PROPOSITION}}
\newtheorem{DEFINITION}[THEOREM]{Definition}
\newenvironment{definition}{\begin{DEFINITION} \hspace{-.85em} {\bf
:} \rm}%
                            {\end{DEFINITION}}
\newtheorem{EXAMPLE}[THEOREM]{Example}
\newenvironment{example}{\begin{EXAMPLE} \hspace{-.85em} {\bf :}
\rm}%
                            {\end{EXAMPLE}}
\newtheorem{CONJECTURE}[THEOREM]{Conjecture}
\newenvironment{conjecture}{\begin{CONJECTURE} \hspace{-.85em}
{\bf :} \rm}%
                            {\end{CONJECTURE}}
\newtheorem{MAINCONJECTURE}[THEOREM]{Main Conjecture}
\newenvironment{mainconjecture}{\begin{MAINCONJECTURE} \hspace{-.85em}
{\bf :} \rm}%
                            {\end{MAINCONJECTURE}}
\newtheorem{PROBLEM}[THEOREM]{Problem}
\newenvironment{problem}{\begin{PROBLEM} \hspace{-.85em} {\bf :}
\rm}%
                            {\end{PROBLEM}}
\newtheorem{QUESTION}[THEOREM]{Question}
\newenvironment{question}{\begin{QUESTION} \hspace{-.85em} {\bf :}
\rm}%
                            {\end{QUESTION}}
\newtheorem{REMARK}[THEOREM]{Remark}
\newenvironment{remark}{\begin{REMARK} \hspace{-.85em} {\bf :}
\rm}%
                            {\end{REMARK}}

\newcommand{\thm}{\begin{theorem}}
\newcommand{\lem}{\begin{lemma}}
\newcommand{\pro}{\begin{proposition}}
\newcommand{\dfn}{\begin{definition}}
\newcommand{\rem}{\begin{remark}}
\newcommand{\xam}{\begin{example}}
\newcommand{\cnj}{\begin{conjecture}}
\newcommand{\mcnj}{\begin{mainconjecture}}
\newcommand{\prb}{\begin{problem}}
\newcommand{\que}{\begin{question}}
\newcommand{\cor}{\begin{corollary}}
\newcommand{\prf}{\noindent{\bf Proof:} }
\newcommand{\ethm}{\end{theorem}}
\newcommand{\elem}{\end{lemma}}
\newcommand{\epro}{\end{proposition}}
\newcommand{\edfn}{\bbox\end{definition}}
\newcommand{\erem}{\bbox\end{remark}}
\newcommand{\exam}{\bbox\end{example}}
\newcommand{\ecnj}{\bbox\end{conjecture}}
\newcommand{\emcnj}{\bbox\end{mainconjecture}}
\newcommand{\eprb}{\bbox\end{problem}}
\newcommand{\eque}{\bbox\end{question}}
\newcommand{\ecor}{\end{corollary}}
\newcommand{\eprf}{\bbox}
\newcommand{\beqn}{\begin{equation}}
\newcommand{\eeqn}{\end{equation}}
\newcommand{\wbox}{\mbox{$\sqcap$\llap{$\sqcup$}}}
\newcommand{\bbox}{\vrule height7pt width4pt depth1pt}
\newcommand{\qed}{\bbox}
\def\sup{^}

\def\H{\{0,1\}^n}
\def\B{\{0,1\}}

\def\S{S(n,w)}

\def\n{\lfloor \frac n2 \rfloor}

\def \E{{\cal E}}
\def \Ex{\mathbb E}
\def \R{\mathbb R}
\def \T{\mathbb T}
\def \Z{\mathbb Z}
\def \F{\mathbb F}

\def \S{{\mathbb S}^{n-1}}
\def \D{{\cal D}}
\def \C{{\cal C}}
\def \H{{\cal H}}
\def \I{{\cal I}}
\def \grad{\mbox{grad}}
\def \Grad{\mbox{Grad}}
\def \div{\mbox{div}}

\def\<{\left<}
\def\>{\right>}
\def \({\left(}
\def \){\right)}
\def \e{\epsilon}
\def \r{\rfloor}

\def \1{{\bf 1}}

\def\Tp{Tchebyshef polynomial}
\def\Tps{TchebysDeto be the maximafine $A(n,d)$ l size of a code with distance $d$hef polynomials}
\newcommand{\rarrow}{\rightarrow}

\newcommand{\larrow}{\leftarrow}

\overfullrule=0pt
\def\setof#1{\lbrace #1 \rbrace}

\begin{abstract}
A {\it design} is a finite set of points in a space on which every "simple" functions averages to its global mean. Illustrative examples of simple functions are low-degree polynomials on the Euclidean sphere or on the Hamming cube.

We prove lower bounds on designs in spaces with a large group of symmetries. These spaces include globally symmetric Riemannian spaces (of any rank) and commutative association schemes with $1$-transitive group of symmetries.

Our bounds are, in general, implicit, relying on estimates on the spectral behavior of certain symmetry-invariant linear operators.
They reduce to the first linear programming bound for designs in globally symmetric Riemannian spaces of rank-$1$ or in distance regular graphs. The proofs are different though, coming from viewpoint of abstract harmonic analysis in symmetric spaces. As a dividend we obtain the  following geometric fact: a design is large because a union of "spherical caps" around its points "covers" the whole space.
\end{abstract}

\section{Introduction}
In the following $M$ is either a compact connected $C^{\infty}$ Riemannian manifold with no boundary (\cite{Chavel1, Helgason}) or a commutative association scheme (\cite{Bannai-Ito, Delsarte}).

A {\it design} (\cite{CS, Delsarte, DL, Dunkl, Levenshtein, Levenshtein1}) is a finite subset of $M$ which is orthogonal to the space of non-constant "simple" functions on $M$.

To define the notions properly, we look at the {\it Laplacian} $\Delta$ on $M$. If $M$ is a Riemannian space, this is the usual Laplacian operator. If $M$ is a graph with degree $d$ and adjacency matrix $A$, we take $\Delta = d \cdot I - A$ to be the usual graph Laplacian. For an association scheme $M$, we pick any relation $R$ in the scheme and take $\Delta$ to be the Laplacian of the associated graph.

In all cases, it is known that the linear operator so defined has nonnegative eigenvalues $0 = \lambda_0 < \lambda_1 < ...$ and all these eigenvalues have finite multiplicity.

\dfn
\label{dfn:design}
A finite subset $\D \subseteq M$ is a {\it design of strength $t$} if for any eigenfunction $\phi$ of the Laplacian which belongs to one of the eigenvalues $0 < \lambda < t$ holds
$$
\sum_{x \in \D} \phi(x) = 0
$$
\edfn

In other words, for a function $f$ in the vector space $V_t$ spanned by the Laplacian eigenfunctions associated with eigenvalues smaller than $t$ holds
$$
\int_M f = \frac{1}{|\D} \sum_{x \in \D} f(x)
$$
This definition is due to \cite{Delsarte} in the context of association schemes and to \cite{Dunkl} in the context of non-discrete compact homogeneous spaces. We observe that in some interesting cases, such as the Euclidean sphere and the Hamming cube, the space $V_t$ is the space of multivariate polynomials of degree at most $d = d(t)$.

We will show a lower bound on the cardinality of a design of strength $t$. To state this bound, we need a notion of an eigenvalue of a subset $\Omega \subseteq M$.

For a Riemannian manifold $M$, $\Omega$ will be a {\it normal domain}, and its eigenvalue $\lambda(\Omega)$ is the minimal eigenvalue of the Laplacian restricted to functions supported on the domain and vanishing on its boundary. In other words, $\lambda(\Omega)$ is the minimal eigenvalue of a function satisfying the Dirichlet boundary conditions on $\Omega$. This value is also referred to as the {\it fundamental tone} of $\Omega$.

For a graph $M$, $\Omega$ could be any subset and $\lambda(\Omega)$ is the minimal eigenvalue of the Laplacian restricted to functions supported on $\Omega$.

Now we can state our bound, in a somewhat vague form. We use $V(\cdot)$ to denote the Riemannian measure if $M$ is a manifold, or a counting measure if $M$ is a graph.

\thm
\label{thm:main-main}
Let $\D$ be a design of strength $t$ on $M$. Let $\Omega$ be a subset of $M$ with eigenvalue $\lambda$. Then, assuming $M$ and $\Omega$ are \verb"sufficiently symmetric" (to be explained later) , and $\lambda < t$, we have
$$
|\D| \ge \frac{t - \lambda}{t} \cdot \frac{V(M)}{V(\Omega)}
$$
In fact, a union of isomorphic copies of $\Omega$ taken around each point of $\D$ essentially covers $M$ (up to a $\frac{t - \lambda}{t}$-factor).

\noindent The precise conditions on $M$ and $\Omega$ are given in Theorems~\ref{thm:main-manifold}~and~\ref{thm:main-graph} below.
\ethm

We remark that the connection between eigenvalues of subsets of $M$ and bounds on designs in $M$ was established in \cite{Friedman-Tillich}, where an implicit version of theorem~\ref{thm:main-main} is proved for the case of the Hamming cube. The usefulness of harmonic analysis in this context was pointed out in \cite{NS1}. Theorem~\ref{thm:main-main} can be viewed as an extension of Proposition~1.3 in \cite{NS2}.

\subsection{Discussion}
To apply the theorem, we need nice subsets $\Omega \subseteq M$, for which we can upper bound the eigenvalue $\lambda(\Omega)$ by a function of the measure $V(\Omega)$.

For manifolds, the Courant nodal domain theorem shows, in particular, that if $\phi$ is an eigenfunction of the Laplacian with eigenvalue $\lambda$, and if $\Omega$ is a nodal domain of $\phi$, then $\lambda(\Omega) = \lambda$. If $M$ is a graph, this claim does not hold, but we only need a one-sided bound $\lambda(\Omega) \le \lambda$ and this is directly verifiable.

Fixing a distinguished point $o \in M$, we may consider 'symmetric' eigenfunctions of the Laplacian, invariant under a subgroup of isometries fixing $o$, and we may take $\Omega$ to be the nodal domain containing $o$. For rank-$1$ symmetric spaces or $1$-transitive distance regular graphs, $\lambda(\Omega)$ and $V(\Omega)$ can be related by the behavior of the first roots of the corresponding family of orthogonal polynomials. Plugging this into Theorem~\ref{thm:main-main}, we recover the first linear programming bound on $1$-transitive distance regular graphs \cite{Levenshtein} and on rank-$1$ symmetric spaces (\cite{Yudin2} for the sphere, and  \cite{Lyubich} for projective spaces).

Another example is the $n$-dimensional flat torus $\T = \R^n/\Z^n$. Let $\Omega$ be a periodization of a Euclidean ball in $\T$. The eigenvalue of $\Omega$ is the same as for a ball in Euclidean space, and is known explicitly. Using Theorem~\ref{thm:main-main}, we recover the bound of \cite{Yudin1} for the torus. In this argument, the lattice $\Z^n$ can be replaced by any lattice $\Lambda \subseteq \R^n$.

For general symmetric spaces, things are more complicated. However, some general inequalities providing an upper bound on $\lambda(\Omega)$ as a function of the geometry of $\Omega$ and $M$ are known, and, in principle, can be used in Theorem~\ref{thm:main-main} (\cite{Buser, Chavel2}).

For general association schemes, things seem to be even more complicated. An essential obstacle is that in the absence of the triangle inequality (provided by the metric structure in distance regular graphs) the action of the Laplacian could drastically increase the support of a function. This leads to difficulties in controlling $\lambda(\Omega)$ as a function of $V(\Omega)$.

\noindent {\bf Covering}: Let $B(r,x)$ denote a metric ball of radius $r$ around $x$ in $M$. Theorem~\ref{thm:main-main} implies that for any $t$, there is a radius $r = r(t)$, such that $M \approx \cup_{x \in \D} B(r,x)$.\footnote{Concentration of measure in the examples below, except (possibly) the Hamming sphere, shows a slight perturbation of the covering provided by the theorem covers the entire space, up to its negligible fraction.}  The function $r\(\cdot\)$ depends on the space, and satisfies $\lambda(B(r(t),x)) < t$. In particular,
\begin{itemize}
\item
Union of spherical caps of radius $r(t)$ around a {\it spherical design} of strength $t$ essentially covers $\S$.
\item
Union of Hamming balls of radius $r(t)$ around a {\it $t$-wise independent set} essentially covers $\{0,1\}^n$.\footnote{This special case was proved in \cite{NS2}.}
\item
Union of Johnson metric balls of radius $r(t)$ around a {\it combinatorial $t$-design} essentially covers the Hamming sphere.
\item
Let $\Lambda \subseteq \R^n$ be a lattice with shortest vector of length $t$. Let $\Lambda^{\ast}$ be the dual lattice. Then {\it any} point in the torus $\T = \R^n/\Lambda^{\ast}$ is a design of strength $t$, since the condition in Definition~\ref{dfn:design} is trivially satisfied. This means that a (periodization of) Euclidean ball $B$ of radius $r(t)$ essentially covers $\T$. Therefore, up to a negligible multiplicative factor, $V(\T) = |\Lambda^{\ast}| \le |B|$, that is $|\Lambda| \ge \frac{1}{|B|}$. Here $|\cdot|$ denotes Euclidean volume. This establishes a lower bound on the covolume of a lattice with shortest vector of length $t$, and recovers a bound of \cite{Levenshtein2} on lattice sphere packing in $\R^n$.

This also implies that a lattice dual to a good sphere-packing lattice is a good {\it lattice quantizer} \cite{CS, Zador}. If $\Lambda$ lies on the best known upper bound for sphere packing in $\R^n$ (\cite{KL}), then, for large $n$, $G\(\Lambda^{\ast}\) \approx 0.062$, close to the optimal value of $0.0586$ (which would have been attained by lattices lying on the bound of \cite{Levenshtein2}).

\end{itemize}

\section{Riemannian manifolds}

\subsection{Preliminaries}
In the following $M$ is a compact connected $C^{\infty}$ Riemannian manifold with no boundary. Let $V$ denote the Riemannian measure on $M$.

We will deal with some function spaces on $M$: the space $L^2(M)$ of all $L^2$-integrable functions on $M$ with the inner product $\<f,g\> = \int_M f g~d V$, and the space $C^k(M)$ of $k$-wise differentiable functions.

We also consider the space ${\cal L}^2(M)$ of all $L^2$-integrable vector fields on $M$ with the inner product $\<X,Y\> = \int_M \<X,Y\> ~dV$, and differentiable fields ${\C}^k(M)$.

There are several naturally defined linear differential operators on differentiable functions and vector fields on $M$. We won't need to give precise definitions of these operators, but rather some of their properties. In the next subsection we collect the needed information about these operators and their properties.

\noindent{\bf Some differential operators on Riemannian manifolds}
\begin{enumerate}
\item {\it The gradient}. For $f \in C^k(M)$, the gradient $\grad(f)$ is a ${\C}^{k-1}$ vector field on $M$.\\
\item {\it The divergence}. For $X \in {\C}^k(M)$, the divergence $\div(X)$ is a $C^{k-1}$ function on $M$.\\
\item {\it The Laplacian}. For $f \in C^k(M)$, the Laplacian $\Delta f$ is a $C^{k-2}$ function on $M$.
$$
\Delta f = \div(\grad(f))
$$
\end{enumerate}

\noindent {\bf Properties}: We have, for a $C^1$ function $f$ and a ${\C}^1$ vector field $X$
$$
\<\grad(f),X\> = - \<f,\div(X)\>
$$

We now recall the definition of the {\it weak derivative}. A vector field $Y \in {\cal L}^2(M)$ is a weak derivative of a function $f \in L^2(M)$ if for any ${\C}^1$ vector field $X$ holds
$$
\<Y,X\> = -\<f, \div(X)\>
$$
We denote $Y = \Grad(f)$.

Let $\H(M)$ be the subspace of $L^2(M)$ consisting of functions with weak derivatives.  On $\H(M)$ we consider a bilinear form (the {\it Dirichlet form}), given by
$$
D[f,g] = \<\Grad(f),\Grad(g)\>
$$
This form is related to the Laplacian by the Green formula. For $f \in C^2(M)$,
$$
D[f,g] = - \<\Delta f, g\>
$$

\noindent{\bf The Laplacian and its eigenfunctions and eigenvalues}

The eigenfunctions of the Laplacian are functions $\phi \in C^2(M)$ satisfying
$$
\Delta \phi = - \lambda \phi
$$
The numbers $\lambda$ on the right hand side of this equality are the {\it eigenvalues} of the Laplacian.

\thm \cite{Chavel1}
\label{thm:Laplacian-eigen}
All the eigenvalues of the Laplacian are nonnegative. The first eigenvalue is $0$, it's of multiplicity $1$, and the eigenfunction associated with it is a constant function. A subspace of eigenfunctions associated with each eigenvalue is finite-dimensional. Eigenspaces associated with different eigenvalues are orthogonal in $L^2(M)$, and $L^2(M)$ is the direct sum of all the eigenspaces. Furthermore, each eigenfunction is in $C^{\infty}(M)$.
\ethm

\subsection{The main theorem for manifolds}
Here we state a version of Theorem~\ref{thm:main-main} for manifolds and begin to prove it.

In the following $\I$ will denote the group of isometries of $M$. We will always assume $\I$ is transitive.

\thm
\label{thm:main-manifold}
Let $\D$ be a design of strength $t$ on $M$. Let $\Omega$ be a normal domain in $M$ with eigenvalue $\lambda$. Then, assuming $M$ has a \verb"nice" group of isometries, and $\Omega$ is invariant under a \verb"nice" subgroup $\I_0$ of isometries, and $\lambda < t$, we have
$$
|\D| \ge \frac{t - \lambda}{t} \cdot \frac{V(M)}{V(\Omega)}
$$
In fact, a union of isometric copies of $\Omega$ taken around each point of $\D$ essentially covers $M$ (up to a $\frac{t - \lambda}{t}$-factor).

Here \verb"niceness" has one of the two following interpretations:
\begin{itemize}
\item
$\I_0$ is a normal subgroup of $\I$, and the factor group ${\I} / {\I}_0$ is abelian.
\item
$M$ is a globally symmetric Riemannian space with a compact group $\I$ of isometries. $\I_0$ is the stabilizer of a point $o \in M$.
\end{itemize}

\ethm

\noindent {\bf Proof:}\\
Let $|\D| = d$. Let $\D = \left\{y_1,\ldots,y_d\right\}$. Let $o \in M$ be a distinguished point. Let $\tau_i$, $i = 1,\ldots,d$, be an isometry of $M$ taking $y_i$ to $o$.

Let $\phi$ be the eigenfunction of $\bar{\Omega}$ associated with eigenvalue $\lambda$. By this, we mean that $\phi$ is the first eigenfunction for the Dirichlet problem on $\Omega$, satisfying $\Delta \phi = -\lambda \phi$ on $\Omega$ and $\phi_{|\partial \Omega} = 0$.

By a version of Theorem~\ref{thm:Laplacian-eigen} for manifolds with boundary, $\phi$ is $C^{\infty}$ on $\bar{\Omega}$. By the Courant nodal domain theorem, $\phi$ is nonnegative on $\Omega$.

Let $f$ be a function on $M$ defined as follows: $f(x) = \phi(x)$ for $x \in \Omega$ and $f(x) = 0$ for $x \not \in \Omega$. Then $f \in \H(M)$ (\cite{Chavel1}). Let $F \in \H(M)$ be defined as follows:
$$
F(x) = \sum_{i=1}^d f\(\tau_i x\)
$$

We will claim, under appropriate conditions, two properties for $F$.
\begin{enumerate}
\item
$$
D[F,F] \le \lambda \<F,F\>
$$
\item
\label{design:orthogonal}
The function $F$ is orthogonal to any eigenfunction $\phi$ of the Laplacian which belongs to one of eigenvalues $0 < \theta < t$.
\end{enumerate}

We observe that these two properties imply the statement of the theorem. In fact, we will show a geometric-flavor statement
\beqn
\label{cover}
|\D| \cdot V(\Omega) \ge V\(\bigcup_{i \in \D} \Omega_i\) \ge V(\mbox{supp}(F)) \ge \frac{t - \lambda}{t} \cdot V(M)
\eeqn
Here $\Omega_i := \tau^{-1}_i\(\Omega\)$. That is, if we take a copy of $\Omega$ around each point of $\D$, we "almost cover" $M$.

We proceed to show (\ref{cover}). By Rayleigh's theorem:
$$
D[F,F] \ge t \cdot \(\<F,F\> - \frac{\<F,1\>^2}{V(M)}\)
$$
This, combined with $D[F,F] \le \lambda \<F,F\>$ implies, via a simple rearrangement
$$
V(M) \cdot \frac{\<F,F\>}{\<F,1\>^2} \le \frac{t}{t - \lambda}
$$
Note that, by Cauchy-Schwarz, and by the definition of $F$, we have
$$
\frac{\<F,F\>}{\<F,1\>^2} \ge \frac{1}{V(\mbox{supp}(F))}
$$
That is,
$$
V(M) \le \frac{t}{t - \lambda} \cdot V(\mbox{supp}(F)),
$$
completing the proof of (\ref{cover}).

It remains to verify when the two claimed properties of $F$ hold.

The first property seems to hold in a rather wide generality.
\lem
\label{lem:up}
In the above assumptions about $M$, and {\it essentially no assumptions} about its group $\I$ of isometries,
$$
D[F,F] \le \lambda \<F,F\>
$$
\elem
\ignore{
\rem
For this lemma we only need the group $\I$ of isometries of $M$ to be transitive and we may take the subgroup ${\I}_0$ stabilizing $\Omega$ to be the trivial subgroup $\I = \{Id\}$.
\erem
}

\prf
We have
$$
D[F,F] = D\left[\sum_{i=1}^d f_i, \sum_{i=1}^d f_i\right] = \sum_{i,j = 1}^d D\left[f_i, f_j\right],
$$
where $f_i(x) := f\(\tau_i x\)$. Now, for $i = j$ we have, by the definition of $f$ and by Green's formula for a manifold with boundary, in our case, $\bar{\Omega}$, recalling that $\phi = 0$ on $\partial \Omega$.
$$
D\left[f_i, f_i\right] = D[f,f] = \int_M \<\Grad(f),\Grad(f)\>~dV =
$$
$$
\int_{\Omega} \<\grad(\phi),\grad(\phi)\>~dV = - \int_{\Omega} \phi \Delta \phi~dV = \lambda \cdot \<\phi, \phi\> = \lambda \cdot \<f,f\>
$$
The third equality is explained by the fact that $\Grad(f)$ is easily checked to equal $\grad(\phi)$ on $\Omega$ and $0$ on its complement.

For $i \not = j$ we have, setting $\Omega_i = \tau^{-1}_i\(\Omega\)$,
$$
D\left[f_i, f_j\right] = \int_{\Omega_i \cap \Omega_j} \<\grad\(\phi_i\), \grad\(\phi_j\) \> dV =
$$
$$
- \int_{\Omega_i \cap \Omega_j} \phi_i \Delta \phi_j dV + \int_{\partial \(\Omega_i \cap \Omega_j\)} \phi_i \cdot \<\phi_j,\nu\> dA = \lambda \cdot \<f_i, f_j\> + \int_{\partial \(\Omega_i \cap \Omega_j\)} \phi_i \cdot \<\phi_j,\nu\> dA
$$
Now, assume that the boundaries $\partial \Omega_i$ and $\partial \Omega_j$ intersect {\it transversally}, that is, by submanifolds of smaller dimension (by perturbing $\Omega_i$ slightly, if needed). Then $\Theta = \Omega_i \cap \Omega_j$ is a normal domain whose boundary is given by, up to lower-dimensional terms,
$$
\partial \Theta = \(\partial \Omega_i\ \cap \Omega_j\) \cup \(\partial \Omega_j\ \cap \Omega_i\)
$$
We have that the outer normal derivative of $\phi_j$ on $\partial \Omega_j$ is non-positive, since $\phi_j$ is nonnegative on $\Omega_j$ and $0$ on $\partial \Omega_j$. Similarly, $\phi_i \ge 0$ on $\Omega_i$ and $0$ on $\partial \Omega_i$. Taken together,
this implies $\int_{\partial \(\Omega_i \cap \Omega_j\)} \phi_i \cdot \<\phi_j,\nu\> dA \le 0$, and consequently
$$
D\left[f_i, f_j\right] \le \lambda \cdot \<f_i, f_j\>,
$$
completing the proof of the lemma.
\eprf

The second property of $F$, which we will refer to as a {\it design-like}, requires more work. In particular, it needs more assumptions on the group of isometries of $M$ and on the domain $\Omega$. We will see some cases in which it holds in the next subsection.

\subsection{The function $F$ is design-like for sufficiently symmetric $M$ and $\Omega$}
In this section we conclude the proof of Theorem~\ref{thm:main-manifold}. The main role is played by the isometry group $\I(M)$ of $M$. We will show that if $\I$ is sufficiently nice, $F$ is design-like.

The next lemma deals with the first case of the theorem.

\lem ($\I$ has an abelian factor)\\
\label{lem:low}
The function $F$ is design-like, provided one of the following two cases holds:
\begin{enumerate}
\item
The group $\I$ of isometries of $M$ is abelian.
\item
More generally, the group $\I$ contains a normal subgroup ${\I}_0$, such that the factor group ${\I} / {\I}_0$ is abelian and $\Omega$ is ${\I}_0$-invariant. In this case, we choose $f$ to be ${\I}_0$-invariant as well.
\end{enumerate}
\elem
\ignore{
\rem
Here we need more from the groups ${\I}$ and ${\I}_0$, but nothing from the domain $\Omega$.
\erem
}

\prf
Let $\phi$ be an eigenfunction that belongs to one of the eigenvalues $0 < \theta < t$. We have
$$
\<F,\phi\> = \sum_{i=1}^d \<f_i,\phi\> = \sum_{i=1}^d \int_{M} f\(\tau_i x\) \phi(x) dV =
$$
$$
\sum_{i=1}^d \int_{M} f(y) \phi\(\tau^{-1}_i y\) dV = \int_M f(y) \cdot \(\sum_{i=1}^d \phi\(\tau^{-1}_i y\)\) dV
$$
Now, consider the first case of the lemma. let $\beta \in {\I}_1$ be an isoperimetry of $M$ taking $y$ to $o$. Then $\beta = \tau^{-1}_i \circ \beta \circ \tau_i:~\tau^{-1}_i y~ \longmapsto ~y_i$.
This means $\left\{\tau^{-1}_i y,~i=1...d\right\}$ is an isometric image of the design $\D$, and therefore a design of strength $t$ as well. Consequently, $\sum_{i=1}^d \phi\(\tau^{-1}_i y\) = 0$ for all $y \in M$, and we are done.

Next, consider the second case of the lemma. Observe that if $\Omega$ is ${\I}_0$-invariant, we can make $f$ to be ${\I}_0$-invariant as well, by averaging it over the compact set ${\I}_0$.\footnote{Note that $f$ and its isometric shifts lie in a finite dimensional eigenspace of the Laplacian on $\bar{\Omega}$.}

We may assume $\phi$ is ${\I}_0$-invariant as well. Indeed, if not, let $\mu$ be the normalized Haar measure on ${\I}_0$ and $\psi(y) = \int_{\alpha \in {\I}_0} \phi(\alpha y) d\mu(\alpha)$. Clearly, $\psi$ is $\I_0$-invariant. In addition, $\psi$ is an eigenfunction of the Laplacian with the same eigenvalue as $\phi$. This is true for any isometry $\alpha \in \I$ and function $\phi_{\alpha}$ defined by $\phi_{\alpha}(y) = \phi(\alpha y)$, and for $\psi$ as an average of $\phi_{\alpha}$ over $\alpha \in {\I}_0$.

Since ${\I}_0$ is a normal subgroup, $\int_{{\I}_0} \phi(\tau^{-1}_i \alpha y) d\mu(\alpha) = \int_{{\I}_0} \phi\(\alpha \tau^{-1}_i y\) d\mu(\alpha) = \psi\(\tau^{-1}_iy\)$. Therefore, since $f$ is invariant, we have
$$
\int_M f(y) \cdot \(\sum_{i=1}^d \phi\(\tau^{-1}_i y\)\) dV = \int_M \(\int_{{\I}_0} f(\alpha y) \cdot \(\sum_{i=1}^d \phi\(\tau^{-1}_i y\)\) d\mu(\alpha)\) dV =
$$
$$
\int_M f(y) \cdot \(\int_{{\I}_0} \(\sum_{i=1}^d \phi\(\tau^{-1}_i \alpha^{-1} y\)\) d\mu(\alpha)\) dV =
\int_M f(y) \cdot \(\sum_{i=1}^d \psi\(\tau^{-1}_i y\)\) dV
$$
Now, since $\I/{\I}_0$ is abelian, and $\psi$ is ${\I}_0$-invariant, we have, for any two isometries $\alpha,\beta \in \I$, that $\psi(\alpha \beta y) = \psi(\beta \alpha y)$. Therefore, taking $\alpha = \tau^{-1}_i$, and $\beta:~y \longmapsto o$, we can proceed as in the first case.
\eprf

\xam
Torus, being an abelian group, and acting on itself, is an example of a manifold satisfying the first case of the lemma.
\exam

The next lemma deals, in particular, with the case of rank-$1$ symmetric spaces, for which the assumptions of the lemma are known to hold (\cite{Helgason}), if we take $\I_0$ to be the stabilizer of a distinguished point $o \in M$.

\lem (Unique invariant eigenfunctions)\\
\label{lem:unique-invariant}
Let $o \in M$ be a distinguished point. Suppose there exists a compact subgroup $\I_0$ that fixes $o$, and such that,
for any eigenvalue $\lambda$ of the Laplacian, there is a \verb"unique" $\I_0$-invariant eigenfunction $\phi_{\lambda}$ satisfying $\phi_{\lambda}(o)) = 1$.

Then, if $\Omega$ is $\I_0$-invariant, $F$ is design-like.
\elem

\prf
In the notation above, and following the same line of reasoning, we have
$$
\<F,\phi\> = \int_M f(y) \cdot \(\int_{\alpha \in {\I}_0} \(\sum_{i=1}^d \phi\(\tau^{-1}_i \alpha^{-1} y\)\) d\mu(\alpha)\) dV
$$
Let $\theta$ be the eigenvalue associated with $\phi$. Let $\psi_i(y) = \int_{\alpha \in {\I}_0} \phi\(\tau^{-1}_i \alpha^{-1} y\) d\mu(\alpha)$. Then, as above, $\psi_i$ is an eigenfunction of the Laplacian with the same eigenvalue $\theta$.

We also claim $\psi_i$ is ${\I}_0$-invariant. Indeed, for any $\alpha_1 \in {\I}_0$, we have, due to left-invariance of $\mu$,
$$
\psi_i(\alpha_1 y) = \int_{\alpha \in {\I}_0} \phi\(\tau^{-1}_i \alpha^{-1} \alpha_1 y\) d\mu(\alpha) = \int_{\beta \in {\I}_0} \phi\(\tau^{-1}_i \beta^{-1} y\) d\mu(\beta) = \psi_i(y)
$$
So, by uniqueness, $\psi_i$ is a scalar multiple of $\phi_{\mu}$, $\psi_i = c_i \cdot \phi_{\mu}$.
This means,
$$
\<F,\phi\> = \(c_1 + \cdots + c_d\) \cdot \int_M f(y) \phi_{\mu}(y) dV
$$
We will now show $\sum_{i=1}^d c_i = 0$, completing the proof.

In fact, consider the value of $\sum_{i=1}^d \psi_i$ at the stable point $o$. We have,
$$
\(\sum_{i=1}^d \psi_i\)(o) = \sum_{i=1}^d \int_{\alpha \in {\I}_0} \phi\(\tau^{-1}_i \alpha^{-1} o\) d\mu(\alpha) =
\sum_{i=1}^d \phi\(\tau^{-1}_i o\) = \sum_{x \in \D} \phi(x) = 0
$$
Therefore,
$$
0 = \(\sum_{i=1}^d \psi_i\)(o) = \(c_1 + \cdots + c_d\) \cdot \phi_{\mu}(o) = c_1 + \cdots + c_d.
$$
\eprf

Given some facts from harmonic analysis on symmetric spaces, the following more general claim holds with essentially the same proof. This claim concludes the proof of Theorem~\ref{thm:main-manifold}.

\lem
\label{lem:symmetric spaces}
Let $M$ be a {\it globally symmetric Riemannian space} with a compact group $\I$ of isometries. Let $\I_0$ be the stabilizer of a point $o \in M$. Then, if the domain $\Omega$ is $\I_0$ invariant, $F$ is design-like.
\elem

\prf (Lemma~\ref{lem:symmetric spaces})
The proof of the lemma proceeds very similarly to the proof of Lemma~\ref{lem:unique-invariant}. We will show that, for any $y \in M$,
\beqn
\label{ident-zero}
\sum_{i=1}^d \psi_i(y) := \sum_{i=1}^d \int_{\alpha \in {\I}_0} \phi\(\tau^{-1}_i \alpha^{-1} y\) d\mu(\alpha) = 0,
\eeqn
for any eigenfunction $\phi$ of the Laplacian with eigenvalue $0 < \theta < t$.

For this, we will need two facts about the algebra ${\bf D}(M)$ of invariant differential operators on symmetric spaces. First, this algebra is commutative, and the joint eigenfunctions of ${\bf}D(M)$ span $L^2(M)$ (\cite{Helgason}, Ch. 10. ex.3). In addition, an $\I_0$-invariant joint eigenfunction $\sigma$ is determined by its value in $o$ (which is non-zero, if $\sigma \not = 0$), and by its eigenvalues on ${\bf D}(M)$ \cite{Helgason}.

Since the function $\phi$ belongs to a finite-dimensional eigenspace $V_{\theta}$ of the Laplacian, and since eigenfunctions associated to distinct eigenvalues are orthogonal, we can write $\phi = \sum_{j=1}^k t_j \sigma_j$ as a linear combination of joint eigenfunctions of ${\bf D}(M)$ with $\sigma_j \in V_{\theta}$. Therefore, it suffices to prove (\ref{ident-zero}) for one joint eigenfunction $\phi$.

Similarly to the above, if for $D \in {\bf D}(M)$ holds $D \phi = \lambda_{\phi} \phi$, also $D \psi_i = \lambda_{\phi} \psi_i$. Therefore, by uniqueness of joint eigenfunctions with the same set of eigenvalues, $\psi_i$ are constant multiples of the $\I_0$-invariant joint eigenfunction $\sigma$, which belongs to the same eigenvalues as $\phi$.

Now, we can conclude the proof similarly to that of Lemma~\ref{lem:unique-invariant}.
\eprf

\section{Graphs}
Let $M = (V,E)$ be a finite graph with $N$ vertices. We assume $M$ is $k$-regular and, moreover, has a $1$-transitive group $\I$ of symmetries.

We define the Laplacian of $M$ to be the linear operator acting on real-valued functions on $M$ in the following manner:
$$
\Delta f(x) =  kf(x) - \sum_{y\sim x} f(y)
$$
We also define the Dirichlet form $D[f,f]$ by setting
$$
D[f,f] = \<f,\Delta f\> = \frac12 \sum_{x\sim y} (f(x) - f(y))^2
$$
Let $0 \le \lambda_1 \le \lambda_2 \le \ldots \lambda_N$ be the eigenvalues of the Laplacian, and let $\phi_1,\phi_2,\ldots \phi_N$ be corresponding eigenfunctions.

We define a design $\D$ of strength $t$ in $M$ exactly as above. We take $d = |\D|$, and set $\D = \left\{y_1,\ldots,y_d\right\}$.

Let $\Omega$ be a subset of $V$. We define the first eigenvalue $\lambda(\Omega)$ to be
$$
\lambda(\Omega) = \min\left\{\frac{D[f,f]}{\<f,f\>}:~f: \Omega \rarrow \R;~f \not = 0\right\}
$$
We remark that in the discrete case we do not require $f$ to vanish on the vertex boundary of $\Omega$ (it "automatically" vanishes on the {\it outer} vertex boundary).

Similarly to the above, we have the following theorem
\thm
\label{thm:main-graph}
Let $\D$ be a design of strength $t$ on $M$. Let $\Omega$ be a subset of $M$ with eigenvalue $\lambda$. Then, assuming $M$ has a \verb"nice" group of symmetries, and $\Omega$ is invariant under a \verb"nice" subgroup $\I_0$ of isometries (and occasionally is even \verb"extra nice"), and $\lambda < t$, we have
$$
|\D| \ge \frac{t - \lambda}{t} \cdot \frac{V(M)}{V(\Omega)}
$$
In fact, a union of isometric copies of $\Omega$ taken around each point of $\D$ essentially covers $M$ (up to a $\frac{t - \lambda}{t}$-factor).

Here \verb"niceness" has one of the following interpretations.
\begin{itemize}
\item
$\I_0$ is a normal subgroup of $\I$, and the factor group ${\I} / {\I}_0$ is abelian.
\item
$M$ is a distance transitive graph with a group of isometries $\I$. $\I_0$ is the stabilizer of a point $o \in M$.
\item
$M$ is a  a distance regular graph with a group of isometries $\I$. $\I_0$ is the stabilizer of a point $o \in M$. In this case we also assume $\Omega$ to be {\it spherical}, that is a union of metric spheres centered at $o$. \footnote{Note that this contains the previous case of distance transitive graphs.}
\item
$M$ is a commutative association scheme with a group of symmetries $\I$. $\I_0$ is the stabilizer of a point $o \in M$. In this case we also assume $\Omega$ to be {\it spherical}, that is a union of adjacency classes of $o$. \footnote{Note that this contains the previous case of distance regular graphs.}
\end{itemize}
\ethm

The proof of the first two cases proceeds almost exactly as above. Let us indicate some minor distinctions and how to deal with them.

First, consider the proof of Lemma~\ref{lem:up}. Here is the argument in the graph case.

Let $\phi$ be the eigenfunction of $\Omega$ corresponding to the eigenvalue $\lambda$. Then $\phi \ge 0$ on $\Omega$ and $\Delta \phi = \lambda \phi$ on $\Omega$. Let $f$ be the extension of $\phi$ onto $G$, with $f$ vanishing outside $\Omega$. Note that this means $f \ge 0$, $\Delta f \le \lambda f$. Now, $F = \sum_{i=1}^d f_i$ be defined as above. Then
$$
D[F,F] = \sum_{i,j} D\left[f_i,f_j\right] = \sum_{i,j} \<\Delta f_i,f_j\> \le \lambda \cdot \sum_{i,j} \<f_i,f_j\> = \lambda \<F,F\>
$$

Second, we remark that conditions of Lemma~\ref{lem:unique-invariant} are satisfied if the group $\I$ of isometries is $2$-transitive. Choose any point $o \in M$, and let $\I_0$ be the stabilizer of $o$. Let the eigenvalue $\lambda$ and a corresponding invariant eigenfunction $L$ be given. Then, we claim $L$ is determined by $l_0 := L\(x_0\)$. Indeed, the value $l_1$ of $L$ at the neighbors $y_1,\ldots,y_k$ of $o$ is determined by the equation $k l_0 - k l_1 = \lambda l_0$, the value $l_2$ of $L$ at the points at distance $2$ from $o$ is determined by $l_0, l_1$ and so on.

The case of distance regular graphs and general association schemes requires more work. We discuss it in the next subsection.

\subsection{Association schemes}
In this subsection we deal directly with general association schemes. Distance regular graphs are a special case, singled out in the statement of the theorem for its intrinsic interest.

We start with interpreting the notion of the algebra of invariant differential operators on $M$ and of spherical and zonal spherical functions on $M$ if $M$ is an association scheme. It is conveniently done using notions from the theory of association schemes.

Let $M$ be an association scheme with $m+1$ classes \cite{Bannai-Ito, Delsarte}.
Let $I = A_0,A_1...,A_m$ be the adjacency matrices of the scheme, and $1 = n_0, k = n_1,...,n_m$ the valencies (degrees) of the graphs defined by $A_i$. Let $E_0,...,E_m$ be the projection matrices of the scheme, with $V_i$ the image of $E_i$ and let $m_i = \mbox{rank}\(E_i\) = \mbox{dim}\(V_i\)$ be the multiplicities of the scheme.

A one-to-one transformation $\tau:~M\rarrow M$ is a {\it symmetry} of the scheme if $x\sim_i y \Leftrightarrow \tau x \sim_i \tau y$ for all $x, y \in M$ and all $i = 0,...,m$. The symmetries of the scheme form a group $\I$ and we assume $\I$ is $1$-transitive on $M$.

We take the algebra ${\bf D}(M)$ of {\it invariant differential operators} on $M$ to be generated by
$\Delta_i$, $i = 1,...,m$ with $\(\Delta_i f\)(x) = n_i f(x) - \sum_{y \sim_i x} f(y)$. In particular, $\Delta_1$ is the discrete Laplacian $\Delta$ defined above. This is a commutative algebra of symmetric matrices.

Fix a distinguished point $o \in M$ and call a function $f$ {\it spherical} is if it is constant on metric spheres around $o$. The spherical functions are spanned by $A_i \delta_o$, where $\delta_o$ is the delta-function at $o$. Another basis for the spherical functions are the {\it zonal spherical functions}
$$
\phi_i = \frac{|M|}{m_i} \cdot E_i \delta_o
$$
By definition, $\phi_i$ is a joint eigenfunction of all the operators in ${\bf D}(M)$. These functions are in fact spherical,  span all spherical functions, and are determined by their eigenvalues on ${\bf D}(M)$. All this is due to the fact that $A_0,...,A_m$ and $E_0,...,E_m$ form two bases of the Bose-Mesner algebra of the association scheme. (Thus the projection matrices both span and are spanned by the adjacency matrices.) Note that they are normalized to be $1$ at $o$.

Consider the orthogonal projection from the space of all real-valued functions on $M$ to the subspace of spherical functions. For a function $f:M~\rarrow \R$, this projection is given by
$$
{\cal S}f = \sum_{i=0}^m \frac{\<A_i \delta_o,f\>}{n_i} \cdot A_i \delta_o
$$
A key property is that any joint eigenfunction projects into a scalar multiple of a zonal spherical function.
\lem
\label{lem:eigen-projection}
Let $\phi \in V_i$. Then
$$
{\cal S} \phi = \phi(o) \cdot \phi_i
$$
\elem
\prf
Let $A_k = \sum_{j=0}^m p_{kj} E_j$, $k=0,...,m$. Then
$$
{\cal S} \phi = \sum_{k=0}^m \frac{\<A_k \delta_o,\phi\>}{n_k} \cdot A_k \delta_o = \sum_{k=0}^m \frac{\<A_k \phi,\delta_o\>}{n_k} \cdot A_k \delta_o = \phi(o) \cdot \sum_{k=0}^m \frac{p_{ki}}{n_k} \cdot A_k \delta_o = \phi(o) \cdot {\cal S}\phi_i = \phi(o) \cdot \phi_i
$$

\eprf

\cor
\label{cor:proj-commutes-BM}
The projection ${\cal S}$ commutes with any operator $A \in {\bf D}(M)$. In other words, ${\cal S}$ commutes with all the operators in the Bose-Mesner algebra.
\ecor
\prf
It suffices to show $A_i {\cal S} = {\cal S} A_i$. Let $f:~M \rarrow \R$ and expand $f = \sum_{j=0}^m a_j f_j$ with $f_j \in V_j$. Then
$$
A_i {\cal S} f = \sum_{j=0}^m a_j f_j(o) A_i \phi_j = \sum_{j=0}^m a_j f_j(o) p_{ij} \phi_j
$$
and
$$
{\cal S} A_i f = {\cal S} \(\sum_{j=0}^m a_j A_i f_j\) = {\cal S} \(\sum_{j=0}^m a_j p_{ij} f_j \) = \sum_{j=0}^m a_j f_j(o) p_{ij} \phi_j
$$
\eprf

The proof of the remaining case of Theorem~\ref{thm:main-graph} can now proceed very similarly to the proof of Lemma~\ref{lem:symmetric spaces}. In the following we refer, both in notation and in argument, to that proof.

First, by Corollary~\ref{cor:proj-commutes-BM}, we may assume the function $f$ to be spherical, since $\Omega$ is. We also may assume the eigenfunction $\phi$ which belongs to eigenvalue $\theta$ of the Laplacian to be a zonal function $\phi \in V_k$. Then, as above, using Lemma~\ref{lem:eigen-projection} and the fact that ${\cal S}$ is an orthogonal projection
$$
\<F, \phi\> = \<f(y), \sum_{i=1}^d \phi\(\tau_i^{-1} y\)\> = \<{\cal S} f(y), \sum_{i=1}^d \phi\(\tau_i^{-1} y\)\> = \<f(y), {\cal S}\(\sum_{i=1}^d \phi\(\tau_i^{-1} y\)\)\> =
$$
$$
\<f, \phi_k\> \cdot \sum_{i=1}^d \phi\(\tau_i^{-1}o\) = \<f, \phi_k\> \cdot \sum_{x \in \D} \phi\(x\) = 0
$$

\end{document}